\documentclass[12pt]{article}
\usepackage{amsmath, amssymb,amsthm, latexsym, times, graphicx}

\def\eps{\varepsilon}
\def\bR{\mathbb{R}}
\def\bC{\mathbb{C}}
\def\bZ{\mathbb{Z}}
\def\bN{\mathbb{N}}
\def\R{\mathbb{R}}

\def\Z{\mathbb{Z}}

\newtheorem{theorem}{Theorem}[section]
\newtheorem{proposition}[theorem]{Proposition}
\newtheorem{lemma}[theorem]{Lemma}

\title{On the Bethe-Sommerfeld conjecture for certain periodic Maxwell operators}
\author{Mariya Vorobets\\ Mathematics Department, Texas A\&M University\\ College Station, TX 77843-3368\\ e-mail mvorobet@math.tamu.edu}
\date{}
\begin{document}

\maketitle

\begin{abstract}
The Bethe-Sommerfeld conjecture states that the spectrum of the stationary Schr\"{o}dinger operator with a periodic potential in dimensions higher than $1$ has only finitely many gaps. After work done by many authors, it has been proven by now in full generality. Another case of a significant interest, due to its importance for the photonic crystal theory, is of a periodic Maxwell operator, where apparently no results of such kind are known. We establish here that in the case of a $2D$ photonic crystal, i.e. of the medium periodic in two variables and homogeneous in the third one, if the dielectric function is separable, the number of spectral gaps of the corresponding Maxwell operator is indeed finite. It is also shown that, as one would expect, when the medium is near to being homogeneous, there are no spectral gaps at all.
\end{abstract}

\section{Introduction}\label{intro}

The Bethe-Sommerfeld conjecture \cite{BetheSomm} states that the spectrum of the stationary Schr\"{o}dinger operator
\begin{equation}\label{E:schroed}
    -\Delta+V(x)
\end{equation}
with a periodic potential $V(x)$ in $\R^n$, when $n\geq 2$, has only finitely many gaps. Starting with \cite{Skrig_DAN2d,DalTrub} and up to \cite{Parnov}, after work done by many authors, it has been proven by now in full generality (see \cite{Skrig_book,Parnov}) for the history and detailed references). In presence of a periodic magnetic potential, the situation becomes much more complex. The corresponding result was proven in $2D$ case in \cite{Mohamed} and \cite{Karp}. The proofs in the latter papers are very technical. In particular, \cite{Mohamed} used microlocal analysis tools of \cite{HelfMoh}. Very recently L.~Parnovski and A.~Sobolev \cite{ParnSob} have settled a much more general case, which allows in particular inclusion of magnetic terms. Another case of a significant interest is of the Maxwell operator in a periodic medium, where apparently no results of such kind are known. The importance of this problem stems from the photonic crystal theory (e.g., \cite{JJWM,FK1,FK2,Kuch}), where existence of spectral gaps is a major issue. We establish here that in the case of a $2D$ photonic crystal, i.e. of the medium periodic in two variables and homogeneous in the third one, if the dielectric function is separable, the number of spectral gaps of the corresponding Maxwell operator is indeed finite. It is also shown that, as one would expect, when the medium is near to being homogeneous, there are no spectral gaps at all.

Let us start with describing the relevant mathematical model. The standard form of the material Maxwell equations is
\begin{equation}\label{E:Maxwell}
    \begin{cases}
    \nabla\cdot {\bf D}=4\pi\rho,\\
\nabla\times {\bf E}=-\frac{1}{c}\frac{\partial {\bf B}}{\partial t},\\
\nabla\cdot {\bf B}=0,\\
\nabla\times {\bf H}=\frac{4\pi}{c}{\bf J}+\frac{1}{c}\frac{\partial {\bf D}}{\partial t}.\\
    \end{cases}
\end{equation}
Here ${\bf E}$ and ${\bf H}$ are electric and magnetic fields, ${\bf D}$ is the electric displacement, ${\bf B}$ - the magnetic induction, $\rho$ - the charge density, $c$ - the speed of light, and ${\bf J}$ - the free current density. The fields ${\bf E}, {\bf H}, {\bf B}, {\bf D}$, and ${\bf J}$ are vector-valued functions from $\bR^3$ (or a subset of $\bR^3$) into $\bR^3$. We will assume absence of free charges and currents, that is $\rho=0$ and ${\bf J}={\bf 0}$.

We are interested in the EM wave propagation in an isotropic dielectric photonic crystal. In this case, the Maxwell equations should be supplemented by the constitutive (or material) equations
\begin{equation}\label{E:constit}
{\bf D}=\eps {\bf E},\, {\bf B}=\mu {\bf H}.
\end{equation}
Here $\eps$ and $\mu$ are scalar time-independent functions called electric permittivity and magnetic permeability, correspondingly. In most photonic crystals considerations it is assumed that the material is nonmagnetic, that is $\mu=1$. We will also assume that the medium is periodic, that is $\eps$ and $\mu$ are periodic with respect to a lattice $\Gamma$ in $\bR^3$. In what follows, we will assume $\Gamma$ to coincide with the integer lattice $\Z^3$.

Under the above assumptions, the Maxwell system reduces to the form
\begin{equation}\label{E:mono}
\left\{ \begin{array}{ll}
\nabla\times {\bf E}=-\frac{1}{c}\frac{\partial {\bf H}}{\partial t},& \nabla\cdot {\bf H}=0,\\
\nabla\times {\bf H}=\frac{1}{c}\eps(x)\frac{\partial {\bf E}}{\partial t},& \nabla\cdot \eps {\bf E}=0.\\
\end{array}\right.
\end{equation}
For mono-chromatic waves of frequency $\omega\in\R$, one has ${\bf \tilde{E}}(x,t)=e^{i\omega t}{\bf E}(x)$, ${\bf\tilde{H}}(x,t)=e^{i\omega t}{\bf H}(x)$, and thus one arrives to the spectral problem
\begin{equation}\label{a1}
\left (
\begin{array}{cc}
0&-\frac{i}{\eps}\nabla^{\times}\\
\frac{i}{\mu}\nabla^{\times}&0\\
\end{array}
\right )
\left (
\begin{array}{c}
{\bf E}\\
{\bf H}\\
\end{array}
\right )=
\frac{\omega}{c}
\left (
\begin{array}{c}
{\bf E}\\
{\bf H}\\
\end{array}
\right )
\end{equation}
on the subspace $\mathcal{S}$ of smooth vector fields
$\left (
\begin{array}{c}
{\bf E}\\
{\bf H}\\
\end{array}
\right )
$
satisfying
\begin{equation}\label{a2}
\nabla\cdot \eps {\bf E}=0,\hskip.2in \nabla\cdot \mu {\bf H}=0.
\end{equation}

The operator
\begin{equation}\label{a3}
M=\left (
\begin{array}{cc}
0&-\frac{i}{\eps}\nabla^{\times}\\
\frac{i}{\mu}\nabla^{\times}&0\\
\end{array}
\right )
\end{equation}
is called the Maxwell operator. We consider $M$ as the operator on the subspace $\mathcal{S}$.

We can extend the operator $M$ to a self-adjoint operator acting on a Hilbert space $\mathcal{H}$. The Hilbert space $\mathcal{H}$ is a closed subspace of $\mathcal{H}_0=L_2(\bR^3,\bC^3,\eps dx)\oplus L_2(\bR^3,\bC^3,\mu dx)$. Namely, $\mathcal{H}$ consists of those vector fields $({\bf E}, {\bf H})\in \mathcal{H}_0$ for which $\nabla\cdot \eps {\bf E}=0$, $\nabla\cdot \mu {\bf H}=0$, where the divergence is understood in the distributional sense. Note that  $\mathcal{H}$ is the closure of $\mathcal{S}\cap \mathcal{H}_0$ in the Hilbert space $\mathcal{H}_0$. The operator $M$ is naturally extended to act on the set
$$D=\{({\bf E}, {\bf H})\in \mathcal{H}\mid\nabla\times {\bf E}\in L_2(\bR^3,\bC^3,\mu dx),\  \nabla\times {\bf H}\in L_2(\bR^3,\bC^3,\eps dx)\}.$$
Here, as before, the differentiations are understood in the distributional sense. Now $M$ is a self-adjoint operator on $\mathcal{H}$ (see \cite{FK2} for more details). In what follows, we only need to know that $\omega/c\in\bR $ is in the spectrum of the Maxwell operator if the system (\ref{a1}) has a nonzero bounded solution $({\bf E},{\bf H})\in\mathcal{S}$.

One of the fields ${\bf E}$ or ${\bf H}$ could be eliminated and the problem can be re-written in terms of another. For instance, one can re-write (\ref{a1}), (\ref{a2})  as the following second order spectral problem:
\begin{equation}\label{E:2ndorder}
 \begin{cases}
\nabla\times\nabla\times {\bf E}=\lambda\eps(x) {\bf E}\\
\nabla\cdot \eps {\bf E}=0,\\
\end{cases}
\end{equation}
where the spectral parameter $\lambda$ is equal to $\left(\frac{\omega}{c}\right)^2$. Note that when $\eps$ and $\mu$ are real-valued functions the spectrum of the Maxwell operator is symmetric with respect to the origin. Therefore $\lambda=\left(\frac{\omega}{c}\right)^2$ is in the spectrum of the generalized spectral problem (\ref{E:2ndorder}) if and only if both $\omega/c$ and $-\omega/c$ are in the spectrum of $M$ (see  \cite{FK2} for the analogous conclusion).

Our principal task is to show that under appropriate conditions on the periodic dielectric function, the spectrum of the problem (\ref{E:2ndorder}), and hence the spectrum of the operator $M$, has only finitely many gaps. While we expect this statement to hold in general, this text is devoted to proving it in a special case when
\begin{equation}\label{E:sep}
\eps (x_1,x_2,x_3)=\eps_1 (x_1)+\eps_2 (x_2).
\end{equation}
We will show that the number of spectral gaps is finite, even if we restrict our consideration to the invariant subspace of the electric fields ${\bf E}$ that are normal to the plane $(x_1,x_2)$ of periodicity and depend on $(x_1,x_2)$ only.

One should notice that, in spite of many similarities, there are some important differences between the spectral problems for Schr\"{o}dinger and Maxwell operators. This difference arises due to the multiplicative rather that additive appearance of the spectral parameter. This, in particular, applies to existence and location of gaps (see, e.g., \cite{Kuch}). It is easy to create gaps at the bottom of the spectrum of a periodic Schr\"{o}dinger operator (for instance, creating a periodic array of well separated potential wells, see \cite{Kuch}). On the other hand, the spectrum of the problem (\ref{E:2ndorder}) always starts at zero, thus preventing a similar gap opening approach.

The paper is structured as follows: the main results (Theorems \ref{main}, \ref{t1}, and \ref{t2}) are stated in Section 2. It is also noticed there that Theorems \ref{t1} and \ref{t2} imply Theorem \ref{main}. In Section 3, the proof of Theorem \ref{t1} is reduced to an auxiliary Proposition \ref{p1.3}. This Proposition, as well as other auxiliary statements, are proven in Sections 4 and 5. Section 6 contains the proof of Theorem \ref{t2}. The final Section 7 is devoted to final remarks and acknowledgments.
\section{Statement of the results}

Under the imposed assumption (\ref{E:sep}), the Maxwell operator (\ref{a3}) admits an invariant subspace $S_0\subset\mathcal{S}$ of fields $({\bf E},{\bf H})=(E_1, E_2, E_3, H_1, H_2, H_3)$ that do not depend on $x_3$. Furthermore, the space $S_0$ is decomposed into the direct sum of two subspaces $S_1\oplus S_2$, where $S_1$ consists of the fields $(E_1,E_2,0,0,0,H)$ and $S_2$ consists of the fields $(0,0,E,H_1,H_2,0)$. In physical terms, $S_1$ consists of the transverse electric (TE) polarized fields, while $S_2$ consists of the transverse magnetic (TM) polarized fields. It is easy to observe that both $S_1$ and $S_2$ are invariant under the operator $M$. To show that the Maxwell operator has finitely many gaps or no gaps at all, it is enough to consider only TM polarized fields. In terms of the spectral problem (\ref{E:2ndorder}), we assume that ${\bf E}=(0,0,E(x_1,x_2))$. Then the problem (\ref{E:2ndorder}) reduces to the 2D scalar spectral problem
\begin{equation}\label{E:scalar}
   -\Delta E=\lambda\eps(x)E.
\end{equation}
Thus the spectrum of problem (\ref{E:2ndorder}) contains the spectrum of the operator $-\dfrac{1}{\eps}\Delta$ considered as a self-adjoint operator on the Hilbert space $L_2(\bR^2, \eps(x_1,x_2)\,dx_1\, dx_2)$.

Our main result is:

\begin{theorem}\label{main}
Let $\eps (x_1,x_2,x_3)=\eps_1 (x_1)+\eps_2 (x_2)$, where $\eps_1$ and $\eps_2$ are $C^2$-smooth, positive, $1$-periodic functions on $\bR$. Then,
\begin{enumerate}
\item The spectrum of the problem (\ref{E:scalar}), and hence of the Maxwell operator $M$, contains a ray and thus has only finitely many gaps.
\item If the functions $\eps_1$ and $\eps_2$ are sufficiently close to constants uniformly on the whole real axis, then the spectrum of (\ref{E:scalar}) coincides with $[0,\infty)$ and has no gaps at all (in this case, the spectrum of $M$ coincides with the whole real axis).
\end{enumerate}
\end{theorem}

The well known Bloch theorem (see, e.g., its most general formulation in \cite[Theorem 4.3.1]{Kuchment}) provides a nice description of the spectrum of elliptic differential operators with periodic coefficients. In our case this theorem can be formulated as follows.

\begin{proposition}(e.g., \cite[Theorem 4.3.1]{Kuchment})\label{l1.1}

Let $\eps\in C^2(\bR^2)$ be a positive function periodic with respect to a lattice $l_1\bZ\oplus l_2\bZ$, i.e.,
$$\eps(x_1+l_1n_1,x_2+l_2n_2)=\eps(x_1,x_2)$$
for all $x_1,x_2\in\bR$ and $n_1,n_2\in\bZ$. Then the following statements are equivalent:

\begin{enumerate}

\item
A number $\lambda\geq 0$ is in the spectrum of the problem (\ref{E:scalar}) (in other, words, it is in the spectrum of the operator $\displaystyle -\frac{1}{\eps}\Delta$).

\item The differential equation
\begin{equation}\label{x1}
-\Delta E=\lambda\eps E
\end{equation}
has a bounded nonzero solution $E$.

\item  The equation (\ref{x1})
has a nonzero Floquet-Bloch solution $E$ that satisfies a cyclic (Floquet) condition
$$E(x_1+l_1n_1,x_2+l_2n_2)=E(x_1,x_2)e^{i(\alpha n_1+\beta n_2)}$$
for some $\alpha,\beta\in\bR$ and all $x_1,x_2\in\bR$, $n_1,n_2\in\bN$.
\end{enumerate}
\end{proposition}

The last of the three statements above is the most convenient for us. We thus prove the following statement, which, according to Proposition \ref{l1.1}, implies the first statement of Theorem \ref{main}.

\begin{theorem}\label{t1}
Let $\eps(x_1,x_2)=\eps_1(x_1)+\eps_2(x_2)$, where $\eps_1$, $\eps_2$ are $C^2$-smooth, strictly positive, $1$-periodic functions on $\bR$. Then there exists $\lambda_0>0$ such that for any $\lambda\geq\lambda_0$ the partial differential equation
$$-\Delta E=\lambda\eps E$$
has a bounded nonzero Floquet-Bloch solution
$$E(x_1,x_2)=E_1(x_1)E_2(x_2),$$
where $E_1$, $E_2$ are such that
\begin{equation}\label{eq1}
E_1(x_1+1)=e^{i\alpha}E_1(x_1),
E_2(x_2+1)=e^{i\beta}E_2(x_2),
\end{equation}
with $\alpha$, $\beta$ $\in \bR$.

Furthermore, $\lambda_0$, depends only on the number
$$
C:=\max\limits_{i=1,2, x\in\R}\{|\eps_i(x)|, |\eps'_i(x)|, |\eps''_i(x)|, |(\eps_i(x))^{-1}|\}.
$$
\end{theorem}

The second statement of Theorem \ref{main} follows if we establish the following result:

\begin{theorem}\label{t2}
Let $\eps\in C(\bR^2)$ be a positive $\Z^2$-periodic function.
Then, for any $\Lambda>0$, there exists $\delta>0$ such that if $|\eps(x)-1|<\delta$, then for any $0\leq \lambda\leq\Lambda$, the partial differential equation
$$-\Delta E=\lambda\eps E$$
has a bounded nonzero solution in $\bR^2$.
\end{theorem}

The particular choice of the period (and thus lattice $\Gamma$) is not important for the proofs and can be made arbitrary by rescaling. For simplicity we will assume, as we have already agreed before, that $\Gamma=\Z^2$, and in particular  ``periodicity'' of a function of one variable, unless specified otherwise, always means ``$1$-periodicity.''

\section{Proof of Theorem \ref{t1}}\label{proof1}

We start with the standard separation of the variables and thus reduction to a one-dimensional problem:

\begin{lemma}\label{l1.2}
Let $\eps_1$, $\eps_2$ be continuous functions on $\bR$ and $\lambda, c\in\bR$. Suppose that $E_1$ is a solution to the differential equation
$$E_1''(x)+\lambda(\eps_1(x)+c)E_1(x)=0$$
and $E_2$ is a solution to the differential equation
$$E_2''(x)+\lambda(\eps_2(x)-c)E_2(x)=0.$$
Then, the function
$E(x_1,x_2)=E_1(x)E_2(x)$
is a solution to the partial differential equation
$$-\Delta E=\lambda\eps E,$$
where $\eps(x_1,x_2)=\eps_1(x_1)+\eps_2(x_2)$.
\end{lemma}
The proof is straightforward.

The proof of Theorem \ref{t1} will be extracted from the following auxiliary result:

\begin{proposition}\label{p1.3}
Let $\eps\in C^2(\bR)$ be a positive $1$-periodic function on $\bR$. Let $d_0>0$ be a constant such that
$$
\max\limits_{x\in\R}\{|\eps(x)|,|(\eps(x))^{-1}|,|\eps'(x)|,|\eps''(x)|\}\leq d_0.
$$
Then there exist positive constants $\lambda_0, d_1, d_2$ that depend only on $d_0$, such that the following property holds:\\
\indent If for some $\lambda\geq\lambda_0$ the equation
$$E''+\lambda\eps E=0$$
does not have any bounded nonzero solutions, then the equation
$$E''+\lambda(\eps+c) E=0$$
has such a solution for any constant $c$ satisfying
$$\frac{d_1}{\lambda}\leq |c|\leq\frac{d_2}{\sqrt{\lambda}}.$$
\end{proposition}

The proof of Proposition \ref{p1.3} will be provided in Section \ref{1b}. Now we are going to show how Theorem \ref{t1} can be derived from this proposition.

\subsection{ Proof of Theorem \ref{t1}}

Let us choose a constant $d_0>0$ such that
$$
\frac{2}{d_0}\leq\eps_i\leq\frac{d_0}{2}, |\eps_i'|\leq d_0, |\eps_i''|\leq d_0 \mbox{ for } i=1,2.
$$
Let also $\lambda_0$, $d_1$, $d_2$ be constants provided by Proposition \ref{p1.3} for this particular $d_0$. We introduce a new constant $$\displaystyle\Lambda_0=\max\left\{\lambda_0, d_0d_1, \frac{2d_1}{d_0}, \left(\frac{2d_1}{d_2}\right)^2\right\}.$$

We will show now that any $\lambda$ greater than $\Lambda_0$ is in the spectrum of (\ref{E:scalar}), which will prove Theorem \ref{t1}.

Let $c_1=0$, $\displaystyle c_2=\frac{d_1}{\lambda}$, and $\displaystyle c_3=-\frac{d_1}{\lambda}$.

We have $\displaystyle |c_j|\leq \frac{d_1}{\Lambda_0}$ for $j=1,2,3$. Since $\Lambda_0\geq d_0d_1$ and $\displaystyle\Lambda_0\geq\frac{2d_1}{d_0}$, we obtain that $\displaystyle |c_j|\leq \frac{d_0}{2}$ and also $\displaystyle |c_j|\leq\frac{1}{d_0}$.
It follows that
$$\frac{1}{d_0}\leq\eps_i+c_j\leq d_0$$
for $i=1,2$ and $j=1,2,3$.
Besides, we clearly have
$$|(\eps_i+c_j)'|=|\eps_i'|\leq d_0$$
and
$$|(\eps_i+c_j)''|=|\eps_i''|\leq d_0.$$
Hence the conclusion of Proposition \ref{p1.3} holds for either of the functions $\eps_i+c_j$.

Let $j,k\in\{1,2,3\}$, $j\ne k$. Then
$$\frac{d_1}{\lambda}\leq |c_j-c_k|\leq \frac{2d_1}{\lambda}.$$
Since $\displaystyle\lambda\geq \Lambda_0\geq \left(\frac{2d_1}{d_2}\right)^2$, we obtain that
$$|c_j-c_k|\leq \frac{2d_1}{\lambda}\leq\frac{2d_1}{\sqrt{\lambda}\sqrt{\Lambda_0}}\leq\frac{d_2}{\sqrt{\lambda}}.$$
Now Proposition \ref{p1.3} implies that the equation
\begin{equation}\label{m1}
E''_1(x)+\lambda(\eps_1(x)+c_j)E_1(x)=0
\end{equation}
does not admit a bounded nonzero solution for at most one value of $j=1,2,3$.

Similarly, the equation
\begin{equation}\label{m2}
E''_2(x)+\lambda(\eps_2(x)-c_j)E_2(x)=0
\end{equation}
does not admit a bounded nonzero solution for at most one value of $j=1,2,3$.
Thus, for at least one $j=1,2,3$ both equations (\ref{m1}) and (\ref{m2}) admit bounded nonzero solutions $E_1$ and $E_2$. Then, according to Lemma \ref{l1.2}, $E(x_1,x_2)=E_1(x_1)E_2(x_2)$ is a solution of the partial differential equation
$$-\Delta E(x_1,x_2)=\lambda(\eps_1(x_1)+\eps_2(x_2))E(x_1,x_2).$$
Since $E$ is clearly bounded and nonzero, Theorem \ref{l1.1} implies that $\lambda$ is in the spectrum of the operator $\displaystyle -\frac{1}{\eps_1+\eps_2}\Delta$. This proves Theorem {\ref{t1}}. \qed

\section{On the spectra of one-dimensional problems}\label{1a}
In order to prove Proposition \ref{p1.3}, we need to conduct an auxiliary study of the spectrum of the one-dimensional differential operator

\begin{equation}\label{E:1D}
-\frac{1}{u}\left(\frac{d^2}{dx^2}-\rho\right),
\end{equation}

\noindent where $u$ and $\rho$ are $l$-periodic functions and $u>0$.
First of all, the Bloch theorem in one-dimensional case implies the following result:

\begin{theorem}(e.g., \cite[Theorem 4.3.1]{Kuchment})\label{l1.4}
Let $u\in C^2(\bR)$ be a positive $l$-periodic function and $\rho$ be piecewise continuous $l$-periodic function. Then the following are equivalent:

\begin{enumerate}

\item
$\lambda$ is in the spectrum of the operator $\displaystyle -\frac{1}{u}\left(\frac{d^2}{dx^2}-\rho\right)$.

\item The differential equation
$$E''(x)+\lambda u(x)E(x)+\rho (x)E(x)=0$$
has a bounded nonzero solution $E$.

\item  For some $\alpha\in\bR$, the cyclic (Floquet) boundary value problem
\begin{equation}\label{aaa1}
\begin{array}{c}
E''(x)+\lambda u(x)E(x)+\rho (x)E(x)=0,\\
E(l)=E(0)e^{i\alpha},
E'(l)=E'(0)e^{i\alpha}\\
\end{array}
\end{equation}
has a nonzero solution.
\end{enumerate}
\end{theorem}

Since the problem is now formulated on a finite interval, the spectrum of (\ref{aaa1}) is no longer continuous. The following lemma is standard:

\begin{lemma}\label{p1.5}
\begin{enumerate}
\item The spectrum of problem (\ref{aaa1}) is discrete and consists of a nondecreasing sequence
$$\lambda_1\leq\lambda_2\leq\lambda_3\leq\ldots$$
such that $\lambda_n\to\infty$.
\item The eigenvalues satisfy the variational principle
$$\lambda_n=\inf\limits_{
\begin{array}{cc}
V\subset H_\alpha^1\\
\dim{V}=n\\
\end{array}}
\sup\limits_{\begin{array}{cc}
f\in V\\
f\ne 0\\
\end{array}}
\frac{f'\cdot f'-(\rho f)\cdot f}{(uf)\cdot f},$$
where $f_1\cdot f_2=\int_0^l f_1(x)\bar{f_2}(x)\ dx$ is the $L^2$-scalar product, $H_\alpha^1$ is a space of functions $f\in H^1[0,l]$ such that  $f(l)=e^{i\alpha}f(0)$, and $V$ is a vector subspace of $H_\alpha^1$.
\item Each eigenvalue $\lambda_n$ depends continuously on the parameter $\alpha\in\bR$, $n=1,2,\ldots$.
    \end{enumerate}
\end{lemma}

Indeed, due to ellipticity and compactness of the interval, the analytic Fredholm theorem (e.g., \cite[Theorem 1.6.16]{Kuchment}) implies that the spectrum either coincides with the whole complex plane, or is discrete. Since obviously large negative values of $\lambda$ are not in the spectrum, the first statement of the lemma follows. The second and third statements are also straightforward.

In view of Lemma \ref{p1.5}, the  range of $\lambda_n$ as a function of $\alpha$ is a closed interval $J_n$, called the \textit{$n$th band of the spectrum}, and the entire spectrum of the operator (\ref{E:1D}) in $L^2(\R)$ is the union of these bands for $n=1,2,...$. The neighboring bands are either adjacent, or else they are separated by a {\it gap}.

Now let us consider problem (\ref{aaa1}) for a fixed value of $\alpha$.
The eigenvalues $\lambda_1, \lambda_2,...$ (and the corresponding eigenfunctions) can be explicitly computed in the special case when $u$ is constant and $\rho=0$. To estimate $\lambda_n$ in a general case, we will reduce the considerations to this special case. The first step here is to apply the Liouville transformation (see \cite{Ol}), which allows one to reduce the problem (\ref{aaa1}) to a similar problem, but now with a constant function $u$ and a different value of the period $l$. This is done in Lemmas \ref{l1.7}--\ref{l1.10}. The function $\rho$ is altered as well, and this is why it has been included in the problem (\ref{aaa1}) in the first place (in the applications we will have $\rho=0$). Then the influence of the function $\rho$ on the spectrum is estimated using the variational principle (see Lemma \ref{l1.8}).

Let $u>0$, $u\in C^2[0,l]$. Let us define a function $$\xi(x)=\int_0^x \sqrt{u(\tau)}\ d\tau.$$ Then $\xi\in C^3[0,l]$ maps the interval $[0,l]$ homeomorphically onto the interval $[0,A]$, where $A=\xi(l)$. We denote by $z$ its inverse function, which is defined on $[0,A]$.

Let
\begin{equation}\label{f1}
\theta=\frac{5[u']^2}{16u^3}-\frac{u''}{4u^2}.
\end{equation}
The function $\theta$ is bounded and its upper bound can be easily estimated in terms of function $u$.

The proof of the following lemma follows by a straightforward calculation.
\begin{lemma}\label{l1.7}
Let $E$ be a function on $[0,l]$ and we introduce a new function $$F(y)=\frac{1}{\sqrt{z'(y)}}E(z(y)).$$
Then,
\begin{enumerate}
\item For any $\lambda\in\bR$, the function $E$ is a solution of the differential equation $$E''(x)+\lambda u(x)E(x)=0, \hskip.2in x\in[0,l]$$
if and only if the function $F$ is a solution to the differential equation $$F''(y)+\lambda F(y)+\theta(z(y))F(y)=0, \hskip.2in y\in[0,A].$$
\item If additionally $u$ satisfies the periodicity conditions $u(l)=u(0), u'(l)=u'(0)$, then the function $F$ satisfies the Floquet conditions $F(A)=F(0)e^{i\alpha},F'(A)=F'(0)e^{i\alpha}$
if and only if $E$ satisfies the similar Floquet conditions $E(l)=E(0)e^{i\alpha}, E'(l)=E'(0)e^{i\alpha}$.
\end{enumerate}
\end{lemma}

Now we are going to compare the spectrum of the problem (\ref{aaa1}) in the case $\rho=0$ with the spectrum of an explicitly solvable problem.

\begin{lemma}\label{l1.8}
Let $u\in C^2[0,l]$ be positive on $[0,l]$ and satisfy (with its first derivative) periodic boundary conditions. Let also $\lambda_1\leq\lambda_2\leq\lambda_3\leq\ldots$ be the eigenvalues of the problem
\begin{equation}\label{bbb1}
\begin{array}{c}
E''(x)+\lambda u(x)E(x)=0,\\
E(l)=E(0)e^{i\alpha}, E'(l)=E'(0)e^{i\alpha}.\\
\end{array}
\end{equation}
Further let $\lambda^*_1\leq\lambda^*_2\leq\lambda^*_3\leq\ldots$
be the eigenvalues of the problem
\begin{equation}\label{bbb2}
\begin{array}{c}
F''(x)+\lambda F(x)=0,\\
F(A)=F(0)e^{i\alpha}, F'(A)=F'(0)e^{i\alpha},\\
\end{array}
\end{equation}
where
\begin{equation}\label{ccc1}
A=\int_0^l \sqrt{u(\tau)}\ d\tau.
\end{equation}
Then $|\lambda_n-\lambda_n^*|\leq \sup|\theta|$ for all $n\geq 1$, where $\theta$ is defined in (\ref{f1}).
\end{lemma}

{\bf Proof } Let, as before, $\xi(x)=\int_0^x\sqrt{u(\tau)}\ d\tau$, $x\in [0,l]$ and let $z$ be its inverse function.

According to Lemma \ref{l1.7}, a function $E$ is the solution to the boundary problem (\ref{bbb1}) if and only if the function $F(y)=\sqrt[4]{u(z(y))}E(z(y))$ is the solution to the boundary problem
\begin{equation}\label{ddd1}
\begin{array}{c}
F''(y)+\lambda F(y)+\rho(y)F(y)=0, \hskip.2in y\in[0,\xi(l)],\\
F(A)=F(0)e^{i\alpha}, F(A)=F(0)e^{i\alpha},\\
\end{array}
\end{equation}
where $\rho(y)=\theta(z(y))$.
Hence we have a one-to-one correspondence between eigenfunctions of the problems (\ref{bbb1}) and (\ref{ddd1}) which is linear and preserves eigenvalues. Therefore, the sequence $\lambda_1\leq\lambda_2\leq\lambda_3\leq\ldots$ is also the spectrum of the problem (\ref{ddd1}).

By the variational principle (Lemma \ref{p1.5}), we have
$$\lambda_n=\inf\limits_{
\begin{array}{cc}
V\subset H_\alpha^1\\
\dim{V}=n\\
\end{array}}
\sup\limits_{\begin{array}{cc}
f\in V\\
f\ne 0\\
\end{array}}
\frac{f'\cdot f'-(\rho f)\cdot f}{f\cdot f},$$
$$\lambda^*_n=\inf\limits_{
\begin{array}{cc}
V\subset H_\alpha^1\\
\dim{V}=n\\
\end{array}}
\sup\limits_{\begin{array}{cc}
f\in V\\
f\ne 0\\
\end{array}}
\frac{f'\cdot f'}{f\cdot f},$$
where $f_1\cdot f_2=\int_0^A f_1(x)\bar{f_2}(x)\ dx$, $H_{\alpha}^1$ is a space of functions $f\in H^1[0,A]$ such that  $f(A)=e^{i\alpha}f(0)$.

Since $\displaystyle\left|\frac{(\rho f)\cdot f}{f \cdot f}\right|\leq \sup|\rho|$ for any nonzero function $f\in H_{\alpha}^1$, it follows that
$$|\lambda_n-\lambda_n^*|\leq \sup|\rho|$$ for all $n\geq 1$.
Clearly, $\sup|\rho|=\sup|\theta|$. This proves the lemma. \qed

Let $\lambda_1\leq\lambda_2\leq\lambda_3\leq\ldots$ be the eigenvalues of the problem (\ref{bbb1}).  We recall that $\lambda_n$ is actually a continuous function of the parameter $\alpha$. Note that the estimate on $\lambda_n$ obtained in Lemma \ref{l1.8} does not depend on $\alpha$. This allows us to estimate the entire band $J_n$, the range of the function $\lambda_n(\alpha)$.

\begin{lemma}\label{l1.10}
Let $C=\sup|\theta|$, where $\theta$ is the function defined by (\ref{f1}), and $A$ be the number defined by (\ref{ccc1}). Then, assuming that $$\frac{\pi^2(n-1)^2}{A^2}+C<\frac{\pi^2n^2}{A^2}-C,$$
one has
$$\left[\frac{\pi^2(n-1)^2}{A^2}+C,\frac{\pi^2n^2}{A^2}-C\right]\subset J_n.$$
Moreover, assuming that $$0<\frac{\pi^2}{A^2}-C,$$
one has the inclusion
$$\left[0,\frac{\pi^2}{A^2}-C\right]\subset J_1$$
\end{lemma}

{\bf Proof } Consider the eigenvalue problem
\begin{equation}
\begin{array}{c}
F''(y)+\lambda F(y)=0,\\
F(A)=F(0)e^{i\alpha},F'(A)=F'(0)e^{i\alpha}.\\
\end{array}
\end{equation}
Its eigenfunctions and eigenvalues are
\begin{equation}
\Lambda_k^*=\left(\frac{\alpha+2\pi k}{A}\right)^2,\,
f_k(y)=\exp\left(\frac{iy(\alpha+2\pi k)}{A}\right), \hskip.2in k\in \bZ.
\end{equation}
Let $\lambda_1^*$, $\lambda_2^*$, ... denote the above eigenvalues arranged in ascending order. Then for $\alpha\in[0,\pi)$ we have
$$\lambda_1^*=\left(\frac{\alpha}{A}\right)^2, \lambda_2^*=\left(\frac{\alpha-2\pi }{A}\right)^2,...,\lambda_{2n-1}^*=\left(\frac{\alpha+2\pi n}{A}\right)^2, \lambda_{2n}^*=\left(\frac{\alpha-2\pi n}{A}\right)^2,...$$
Besides, $\lambda_n(2\pi-\alpha)=\lambda_n(\alpha)$ for all $n\geq 1$ and all $\alpha$.

Let $J_n^*$ be the range of $\lambda_n^*$ as a function of $\alpha$, i.e.
$$
J_n^*=\left[\frac{\pi^2(n-1)^2}{A^2},\frac{\pi^2n^2}{A^2}\right].
$$
In particular, we have $\lambda_n^*(\alpha_1)=\pi^2(n-1)^2 A^{-2}$ and $\lambda_n^*(\alpha_2)=\pi^2 n^2 A^{-2}$ for $\alpha_1=0$, $\alpha_2=\pi$ if $n$ is odd and $\alpha_1=\pi$, $\alpha_2=0$ if $n$ is even.

By Lemma \ref{l1.8}, $|\lambda_n(\alpha_1)-\lambda_n^*(\alpha_1)|\leq C$ and $|\lambda_n(\alpha_2)-\lambda_n^*(\alpha_2)|\leq C$.
Thus the first statement of the lemma follows.
To prove the second statement, it remains to notice  that $0$ always belongs to $J_1$, since $\lambda_1(0)=0$. \qed

\section{Proof of Proposition \ref{p1.3}}\label{1b}

Before proceeding to the proof of Proposition \ref{p1.3}, we need two more technical statements.

Suppose that some $\lambda>0$ lies in a spectral gap of the operator $ -u^{-1}d^2/dx^2$. We would like to know whether it is possible to modify function $u$ slightly, so that $\lambda$ becomes covered by a spectral band. Lemma \ref{l1.10} tells us that $\lambda$ is close to a number $\pi^2n^2A^{-2}$, where  $n$ is an integer and $A$ is defined in (\ref{ccc1}). The next lemma explains how we need to change $A$ to move all such numbers away from $\lambda$. The parameter $b$ in the lemma should be understood as $\pi^2A^{-2}$.

\begin{lemma}\label{l1.11}
Let $b_0$, $b_1>0$, $c>0$. Suppose that
$$|\lambda-bn^2|<c$$
for some $b_0\leq b\leq b_1$ and integer $n\geq 1$. Then $$|\lambda-\tilde{b}m^2|>c$$
for any integer $m\geq 1$, provided that
$$\frac{4b_1c}{\lambda}<|b-\tilde{b}|<\frac{b_0^{3/2}}{12\sqrt{\lambda}}$$
and $\lambda\geq\max\left\{2c,3600c^2 b_0^{-1}\right\}$.
\end{lemma}

{\bf Proof } Assume that $\tilde{b}$ satisfies conditions of the lemma. We show first that $|\lambda-\tilde{b}n^2|>c$.

Indeed, since $\lambda\geq 2c$, the inequality $|\lambda-bn^2|<c$ implies that $\displaystyle \frac{1}{2}<\frac{\lambda}{bn^2}<2$. Then
$$|\lambda-\tilde{b}n^2|\geq |b-\tilde{b}|n^2-|\lambda-bn^2|\geq \frac{4b_1c}{\lambda}n^2-c\geq \frac{4bcn^2}{\lambda}-c>2c-c=c.$$
Now we are going to show that $\lambda-\tilde{b}(n+1)^2<-c$ and $\lambda-\tilde{b}(n-1)^2>c$. This will complete the proof.

We have
$$\lambda-\tilde{b}(n+1)^2=\lambda-b(n+1)^2+(b-\tilde{b})(n+1)^2<c-b(2n+1)+\frac{b_0^{3/2}}{12\sqrt{\lambda}}(n+1)^2\leq c-2bn+\frac{b_0^{3/2}n^2}{3\sqrt{\lambda}}=$$
$$=c-\frac{2bn\sqrt{\lambda}}{\sqrt{\lambda}}+\frac{b_0^{3/2}n^2}{3\sqrt{\lambda}}\leq
c-\sqrt{\frac{b\lambda}{2}}+\frac{\sqrt{b_0}a\sqrt{\lambda}n^2}{3\lambda}\leq
c-\sqrt{\frac{b\lambda}{2}}+\frac{2}{3}\sqrt{b_0\lambda}\leq c-\sqrt{\frac{b_0\lambda}{2}}+\frac{2}{3}\sqrt{b_0\lambda}\leq$$ $$c-\sqrt{b_0\lambda}\left(\frac{2}{3}-\frac{1}{\sqrt{2}}\right).$$
By a hypothesis of the lemma $\sqrt{b_0\lambda}\geq 60c$, hence
$$c-\sqrt{b_0\lambda}\left(\frac{2}{3}-\frac{1}{\sqrt{2}}\right)\leq c-60c\left(\frac{2}{3}-\frac{1}{\sqrt{2}}\right)<-c.$$

It remains to prove that $\lambda-\tilde{b}(n-1)^2>c$.
Since $\lambda\geq 2c$ it is no loss to assume that $n\geq 2$. We have
$$\lambda-\tilde{b}(n-1)^2=\lambda-b(n-1)^2+(b+\tilde{b})(n-1)^2>-c+b(2n-1)-\frac{b_0^{3/2}}{12\sqrt{\lambda}}(n-1)^2.$$
Since $n\geq 2$, we have $(n-1)^2\leq n^2$ and $2n-1\geq n$, then
$$-c+b(2n-1)-\frac{b_0^{3/2}}{12\sqrt{\lambda}}(n-1)^2\geq -c+bn-\frac{b_0^{3/2}}{12\sqrt{\lambda}}n^2\geq -c+\frac{bn\sqrt{\lambda}}{\sqrt{\lambda}}-\frac{\sqrt{b_0}b\sqrt{\lambda}n^2}{12\lambda}\geq $$
$$-c+\sqrt{\frac{b\lambda}{2}}-\frac{1}{6}\sqrt{b_0\lambda}\geq -c+\sqrt{\frac{b_0\lambda}{2}}-\frac{1}{6}\sqrt{b_0\lambda}=-c+\sqrt{b_0\lambda}\left(\frac{1}{\sqrt{2}}-\frac{1}{6}\right)\geq$$
$$ -c+60c\left(\frac{1}{\sqrt{2}}-\frac{1}{6}\right)\geq 30c.$$
This finishes the proof of Lemma \ref{l1.11}. \qed

Let $\eps_0$ be a continuous positive function on $[0,1]$. For any $c>-\inf\eps_0$, let
$$A(c)=\int_0^1 \sqrt{\eps_0(x)+c}\ dx.$$
As we know from Lemma \ref{l1.10}, the quantity $\pi^2A(c)^{-2}$ is closely related to the location of the spectral gaps of the operator $\displaystyle -\eps_0^{-1}d^2/dx^2$. We need to know how it depends on $c$.

\begin{lemma}\label{l1.12}
Suppose that $d_1\leq \eps_0\leq d_2$, where $d_1$ and $d_2$ are positive constants. Then for any $c\in\left [-\frac{d_1}{2}, \frac{d_1}{2}\right ]$ we have
$$\frac{\pi^2}{6d_2^2}|c|\leq \left | \frac{\pi^2}{A^2(c)}-\frac{\pi^2}{A^2(0)}\right |\leq \frac{\pi^2}{2d_1^2}|c|$$
and
$$\frac{2\pi^2}{3d_2}\leq\frac{\pi^2}{A^2(c)}\leq\frac{2\pi^2}{d_1}.$$
\end{lemma}

{\bf Proof } The function $A$ is well defined on $\left [-\frac{d_1}{2}, \frac{d_1}{2}\right ]$ and smooth. We have
$$A'(c)=\frac{d}{dc}\left (\int_0^1 \sqrt{\eps_0(x)+c}\ dx\right )=\int_0^1 \frac{d}{dc}( \sqrt{\eps_0(x)+c})\ dx=\int_0^1 \frac{dx}{2\sqrt{\eps_0(x)+c}}.$$
It follows that $(2\sqrt{d_2+c})^{-1}\leq A'(c)\leq(2\sqrt{d_1+c})^{-1}$. Also, $\sqrt{d_1+c}\leq A(c)\leq\sqrt{d_2+c}$.

Let us introduce a new function $B(c)=\pi^2A^{-2}(c)$.
Then $$B'(c)=-\frac{\pi^2 A'(c)}{2A^3(c)}.$$
Therefore,
$$\frac{\pi^2}{4(d_2+c)^2}\leq |B'(c)|\leq \frac{\pi^2 }{4(d_1+c)^2}.$$
As a consequence,
$$|B(c)-B(0)|\leq \left |\int_0^c\frac{\pi^2}{4(d_1+\tau)^2}\ d\tau\right |\leq\frac{\pi^2 |c|}{4d_1(d_1+c)}.$$
Similarly,
$$|B(c)-B(0)|\geq \frac{\pi^2 |c|}{4d_2(d_2+c)}.$$
In the case $|c|<d_1/2$, we obtain
$$|B(c)-B(0)|\leq \frac{\pi^2 |c|}{4d_1(d_1-\frac{d_1}{2})}=\frac{\pi^2 |c|}{2d_1^2}$$
and
$$|B(c)-B(0)|\geq \frac{\pi^2 |c|}{4d_2(d_2+\frac{d_1}{2})}\geq\frac{\pi^2 |c|}{6d_2^2}.$$
Also, in this case we have
$$\sqrt{d_1/2}\leq\sqrt{d_1+c}\leq A(c)\leq\sqrt{d_2+c}\leq\sqrt{3d_2/2},$$
hence
$$\frac{2\pi^2}{3d_2}\leq B(c)\leq\frac{2\pi^2}{d_1},$$
which proves the lemma. \qed

Now we are ready to prove the central Proposition \ref{p1.3}.

\paragraph{Proof of Proposition \ref{p1.3}}
Recall that the function $\eps$ satisfies $d_0^{-1}
\leq\eps\leq d_0$, $|\eps'|\leq d_0$, $|\eps''|\leq d_0$. Let us introduce constants
$$a_0=\frac{2\pi^2}{3d_0}, \qquad a_1=2\pi^2 d_0, \qquad \Theta=\frac{5}{2}d_0^5+d_0^3,$$ $$d_1=\frac{24}{\pi^2}a_1d_0^2\Theta, \qquad d_2=\frac{a_0^{3/2}}{6\pi^2d_0^2}, \qquad \lambda_0=\max\left\{2\Theta, \frac{3600\Theta^2}{a_0},(2d_0d_2)^2\right\}.$$

Suppose that for some $\lambda\geq \lambda_0$ the differential equation
\begin{equation}\label{f4}
E''+\lambda(\eps+c)E=0
\end{equation}
does not admit a bounded nonzero solution for $c=0$. We shall show that this equation does admit such a solution for any constant $c$ satisfying
$$\frac{d_1}{\lambda}\leq |c|\leq\frac{d_2}{\sqrt{\lambda}}.$$

First we introduce the following functions defined at least for $c>-d_0^{-1}$:
$$A(c)=\int_0^1 \sqrt{\eps(\tau)+c}\ d\tau,$$
and
$$\theta(x,c)=\frac{5[\eps'(x)]^2}{16(\eps(x)+c)^3}-\frac{\eps''(x)}{4(\eps(x)+c)^2}.$$

We shall need some estimates regarding $A$ and $\theta$.

Take any $c$ such that $d_1/\lambda\leq |c|\leq d_2/\sqrt{\lambda}$. Notice that $|c|\leq 1/(2d_0)$ since $\lambda\geq\lambda_0\geq (2d_0 d_2)^2$. Then Lemma \ref{l1.12} implies that
$$\frac{\pi^2}{6d_0^2}|c|\leq \left | \frac{\pi^2}{A^2(c)}-\frac{\pi^2}{A^2(0)}\right |\leq \frac{\pi^2d_0^2}{2}|c|$$
and
$$a_0\leq \frac{\pi^2}{A^2(c)}\leq a_1.$$

Moreover, since $\eps+c\geq \eps-|c|\geq \frac{1}{2d_0}$, $|\eps'|\leq d_0$, and $|\eps''|\leq d_0$, we have
$$|\theta|\leq\left |5\frac{[\eps']^2}{16(\eps+c)^3}\right|+\left |\frac{\eps''}{4(\eps+c)^2}\right | \leq \frac{5d_0^2}{16(2d_0)^{-3}}+\frac{d_0}{4(2d_0)^{-2}}=\Theta.$$

We assumed that the equation (\ref{f4}) does not admit a nonzero bounded solution for $c=0$. By Theorem \ref{l1.4}, $\lambda$ is not an eigenvalue of the operator $-\frac{1}{\eps}\frac{d^2}{dx^2}$. According to Lemma \ref{l1.10}, this implies that
$$\left|\lambda-\dfrac{\pi^2 n^2}{A(0)^2}\right|<\sup\limits_{x\in [0,1]}|\theta(x,0)|\leq\Theta$$ for some integer $n\geq 1$.
By the above we have
$$\left | \frac{\pi^2}{A^2(c)}-\frac{\pi^2}{A^2(0)}\right |\geq \frac{\pi^2}{6d_0^2}|c|\geq \frac{\pi^2}{6d_0^2}\frac{d_1}{\lambda}=\frac{4a_1\Theta}{\lambda}$$
and
$$\left| \frac{\pi^2}{A^2(c)}-\frac{\pi^2}{A^2(0)}\right |\leq \frac{\pi^2d_0^2}{2}|c|\leq\frac{\pi^2d_0^2}{2}\frac{d_2}{\sqrt{\lambda}}=\frac{a_0^{3/2}}{12\sqrt{\lambda}}.$$
Since $\lambda>\max\left\{2\Theta, \frac{3600\Theta^2}{a_0}\right\}$, it follows from Lemma \ref{l1.11} that
$$\left |\lambda-\frac{\pi^2m^2}{A^2(c)}\right |>\Theta$$
for any integer $m\geq 1$.

Since $\sup\limits_{x\in [0,1]} |\theta(x,c)|\leq\Theta$, Lemma \ref{l1.10} implies that $\lambda$ is an eigenvalue of the operator $-\frac{1}{\eps+c}\frac{d^2}{dx^2}$ or, equivalently, the equation (\ref{f4}) has a bounded nonzero solution. \qed

\section{Proof of Theorem \ref{t2}}\label{proof2}

Let $\eps$ be a positive continuous function on $\bR^2$ periodic with respect to the integer lattice. Suppose that $\lambda\geq 0$ belongs to the spectrum of the operator $\displaystyle -\frac{1}{\eps}\Delta$ or, equivalently, the differential equation $$-\Delta E=\lambda\eps E$$
has a bounded nonzero solution in $\bR^2$. According to the Bloch theorem (Theorem \ref{l1.1}), we can choose the bounded solution that satisfies a Floquet condition
\begin{equation}\label{mm1}
E(x_1+l_1n_1,x_2+l_2n_2)=E(x_1,x_2)e^{i(\alpha n_1+\beta n_2)}
\end{equation}
for some $\alpha,\beta\in\bR$ and all $x_1,x_2\in\bR$, $n_1,n_2\in\bN$.
Then the function $E$ is a solution of the following boundary value problem in the unit square with quasiperiodic (or Floquet) boundary conditions:
\begin{equation}\label{l1}
\begin{array}{c}
-\Delta E=\lambda\eps E,\\
E(1,x_2)=e^{i\alpha}E(0,x_2),
\frac{\partial E}{\partial x_1}(1,x_2)=e^{i\alpha}\frac{\partial E}{\partial x_1}(0,x_2),\\
E(x_1,1)=e^{i\beta}E(x_1,0),
\frac{\partial E}{\partial x_2}(x_1,1)=e^{i\beta}\frac{\partial E}{\partial x_2}(x_1,0).\\
\end{array}
\end{equation}
Conversely, any solution of the boundary value problem (\ref{l1}) can be extended to a solution of the equation $-\Delta E=\lambda\eps E$ in the entire plane that satisfies the Floquet condition (\ref{mm1}).

The next standard statement collects the information about the spectrum of the problem (\ref{l1}) that we will need to prove Theorem \ref{t2}.

\begin{proposition}(\cite{Kuchment,RS})\label{p2.1}
\begin{enumerate}
\item The spectrum of the problem (\ref{l1}) is discrete. Its eigenvalues form an nondecreasing sequence
$$0\leq\lambda_1\leq\lambda_2\leq\lambda_3\leq\ldots,\hskip.2in \lambda_n\to\infty.$$
\item Each eigenvalue $\lambda_n$ is a continuous function of $\alpha,\beta\in \bR$.
\item Dependence of the eigenvalue $\lambda_n=\lambda_n(\alpha,\beta;\eps)$ on the function $\eps$ is monotone. Namely, if $\eps\leq\tilde{\eps}$ everywhere in the unit square, then $\lambda_n(\alpha,\beta;\tilde{\eps})\leq\lambda_n(\alpha,\beta;\eps)$.
\end{enumerate}
\end{proposition}
The next lemma provides, also a standard, statement on dependence on the dielectric function $\eps$.
\begin{lemma}\label{l2.3}
For any $\vartheta>0$ and $n$ there exists $\delta>0$ such that
$$|\lambda_n(\alpha,\beta;\tilde{\eps})-\lambda_n(\alpha,\beta;\eps)|<\vartheta$$
for all $\alpha,\beta$ provided that $|\tilde{\eps}-\eps|<\delta$.
\end{lemma}

{\bf Proof. } Clearly,
$$\lambda_n(\alpha,\beta;k\eps)=\frac{1}{k}\lambda_n(\alpha,\beta;\eps)$$
for any $k>0$.
Let us pick $\kappa>0$ such that
$\kappa \lambda_n(\alpha,\beta;\eps)<\vartheta$.
Then
$$\lambda_n(\alpha,\beta;(1+\kappa)\eps)\geq\frac{1}{1+\kappa}\lambda_n(\alpha,\beta;\eps)>(1-\kappa)\lambda_n(\alpha,\beta;\eps)>\lambda_n(\alpha,\beta;\eps)-\vartheta.$$
Also
$$\lambda_n\left(\alpha,\beta;\frac{\eps}{1+\kappa}\right)\leq(1+\kappa)\lambda_n(\alpha,\beta;\eps)\leq\lambda_n(\alpha,\beta;\eps)+\vartheta.$$
If the function $\tilde{\eps}$ satisfies inequalities
$$\frac{\eps}{1+\kappa}\leq\tilde{\eps}\leq(1+\kappa)\eps$$
everywhere in the unit square, then, by Proposition \ref{p2.1},
$$\lambda_n(\alpha,\beta;(1+\kappa)\eps)\leq\lambda_n(\alpha,\beta;\tilde{\eps})\leq\lambda_n\left(\alpha,\beta;\frac{\eps}{1+\kappa}\right),$$
which implies that
$$|\lambda_n(\alpha,\beta;\tilde{\eps})-\lambda_n(\alpha,\beta;\eps)|<\vartheta.$$
Since $$\inf{\eps}=\mu>0,$$ we have
$$(1+\kappa)\eps-\eps=\kappa\eps>\kappa\mu,$$
$$\eps -\frac{\eps}{1+\kappa}=\frac{\kappa\eps}{1+\kappa}\geq\frac{\kappa\mu}{1+\kappa}.$$
Thus, for any function $\tilde{\eps}$ such that $$|\tilde{\eps}-\eps|\leq\frac{\kappa\mu}{1+\kappa}<\kappa\mu,$$
we have $|\lambda_n(\alpha,\beta;\tilde{\eps})-\lambda_n(\alpha,\beta;\eps)|<\vartheta$.
This proves the statement of the lemma. \qed

Let us fix the function $\eps$. For any $n$ let $I_n=I_n(\eps)$ denote the range of $\lambda_n(\alpha,\beta;\eps)$ as a function of $\alpha$ and $\beta$.
It follows from Proposition \ref{p2.1} that $I_n$ is a closed interval. Note that this interval lies in the spectrum of the operator $\displaystyle -\frac{1}{\eps}\Delta$ acting on the entire plane. Furthermore, the spectrum is exactly the union of the intervals $I_1(\eps),I_2(\eps),\ldots$.

The following statement about the spectrum of the Laplace operator is well known (e.g., \cite{Skrig_book}), and can be proven easily, so we skip its proof:
\begin{lemma}\label{l2.4}
In the case $\eps=1$, the intervals $I_1$, $I_2$, $I_3$,... overlap. That is, for any $n$ the intersection $I_n\cap I_{n+1}$ has a nonempty interior.
\end{lemma}

\paragraph{Proof of Theorem \ref{t2}.}
We can now address the proof of the remaining result, Theorem  \ref{t2}, which follows rather immediately from Lemmas \ref{l2.3} and \ref{l2.4}. Indeed:

The spectrum of the operator $\displaystyle -\eps^{-1}\Delta$ is the union of the intervals $I_1(\eps), I_2(\eps),\ldots$ defined earlier in this section. Hence we need to show that for any $\Lambda>0$ this union covers the interval $[0,\Lambda]$ provided that the function $\eps$ is close enough to 1 uniformly.

Let $I_1(1)=[a_1,b_1], I_2(1)=[a_2,b_2],\ldots$. According to Lemma \ref{l2.4}, $a_{n+1}<b_n$.  Given $\Lambda>0$, let us take $N$ such that $b_N>\Lambda$. Let $$\vartheta_0=\frac{1}{2}\min\limits_{1\leq n< N}(b_n-a_{n+1}).$$
By definition, $2\vartheta_0$ is a lower bound on the length of the intersection $I_n\cap I_{n+1}$ for $n=1,2,\ldots,N-1$. Further, let $\vartheta=\min\{\vartheta_0, \Lambda-b_N\}$.

We have $a_n=\lambda_n(\alpha_n,\beta_n;1)$ and $b_n=\lambda_n(\alpha_n',\beta_n';1)$ for some $\alpha_n,\beta_n,\alpha'_n, \beta'_n\in\bR$.
By Lemma \ref{l2.3}, there exists $\delta>0$ such that
$$|\lambda_n(\alpha_n,\beta_n;\eps)-a_n|<\vartheta,$$
$$|\lambda_n(\alpha'_n,\beta'_n;\eps)-b_n|<\vartheta$$
for $n=1,2,\ldots,N$ whenever the function $\eps$ satisfies $\sup|\eps-1|<\delta$. Since the points $\lambda_n(\alpha_n,\beta_n;\eps)$ and $\lambda_n(\alpha'_n,\beta'_n;\eps)$ lie in the interval $I_n(\eps)$, we obtain  $\displaystyle I_n(\eps)\supset [a_n+\vartheta,b_n-\vartheta]$. Moreover, $\displaystyle I_1(\eps)\supset [0,b_1-\vartheta]$ as $\lambda_1(0,0;\eps)=0$ (indeed, the constant function is an eigenfunction of the problem (\ref{l1}) with periodic boundary conditions for $\lambda=0$). By the choice of $\vartheta$, the intervals $I_n(\eps)$ and $I_{n+1}(\eps)$ overlap for $n=1,2,\ldots,N-1$. Besides, the right end of the interval $I_N(\eps)$ lies to the right of the point $\Lambda$.
Thus the intervals $I_1(\eps)$, $I_2(\eps)$,...,$I_N(\eps)$ cover the interval $[0,\Lambda]$ without any gaps. This finishes the proof of Theorem \ref{t2}. \qed

\section*{Final remarks and acknowledgments}

In this paper, we only considered the $E$-polarized modes for electromagnetic waves propagating along the periodicity plane of a $2D$ photonic crystal with a separable dielectric function. There clearly remain several issues to consider. The separability condition is a strong restriction, and thus one would want to avoid it. Besides, the case of fully $3D$ periodic photonic crystals has not been considered. The author plans to address these questions in the future work.

The author is grateful to Peter Kuchment and Yaroslav Vorobets for helpful discussions.


\begin{thebibliography}{99}

\bibitem{BetheSomm}  G. Bethe and A. Sommerfeld, Elektronentheorie der
Metalle, Springer Verlag, Berlin-New York, 1967.

\bibitem{DalTrub}B. E. J. Dahlberg and E. Trubowitz, A remark on two-dimensional periodic potentials, Comment. Math. Helv. {\bf 57} (1982), 130--134.

\bibitem{Eastham}
M. S. P. Eastham, The spectral theory of periodic differential equations, Scottish Acad. Press, Edinburg-London, 1973.

\bibitem{FK1}
A. Figotin and P.Kuchment, Band-gap structure of spectra of periodic dielectric and acoustic media. I: Scalar model, SIAM J. Appl. Math., {\bf 56} (1996), 68--88.


\bibitem{FK2}
A. Figotin and P.Kuchment, Band-gap structure of spectra of periodic dielectric and acoustic media. II: 2D photonic crystals, SIAM J. Appl. Math. {\bf 56} (1996), 1561--1620.

\bibitem{HelfMoh}B.~Helffer and A.~Mohamed,  Asymptotics of the density of states for the Schr�dinger operator with periodic electric potential, Duke Math. J. {\bf 92} (1998), 1--60.

\bibitem{JJWM}
J.~D.~Joannopoulos, S.~G.~Johnson, J.~N.~Winn, and R.~D.~Meade, Photonic Crystals: Molding the Flow of Light, 2nd edition, Princeton Univ. Press, Princeton, 2008.

\bibitem{Karp} Y.~Karpeshina, Spectral Properties of the Periodic Magnetic Schr\"{o}dinger Operator in the High-Energy Region. Two-Dimensional Case, 	Comm. Math. Physics {\bf 251} (2004), No. 3, 473--514.

\bibitem{Kuchment}
P. Kuchment, Floquet theory for partial differential equations. Operator theory advances and applications, vol. 60. Birkhauser, Basel-Boston-Berlin, 1993.

\bibitem{Kuch}
P. Kuchment, The Mathematics of Photonic Crystals, Ch. 7 in "Mathematical Modeling in Optical Science", Gang Bao et al. (Eds.), Frontiers in Applied Mathematics v. 22, SIAM, 2001, 207--272.

\bibitem{Mohamed}A.~Mohamed,  Asymptotics of the density of states for the Schr\"{o}dinger operator with periodic electromagnetic potential, J. Math. Phys. 38 (1997), 4023-4051.

\bibitem{Ol}
F. W. J. Olver, Asymptotics and special functions, Academic Press, NY, 1974.

\bibitem{Parnov} L.~Parnovski, Bethe-Sommerfeld conjecture, Ann. Henri Poincare {\bf 9} (2008), No. 3, 457--508.

\bibitem{ParnSob} L.~Parnovski, A.~V.~Sobolev, Bethe-Sommerfeld conjecture for periodic operators with strong perturbations, Inventiones Math. {\bf 181} (2010), no. 3, 467--540.

\bibitem{RS}
M. Reed and B. Simon, Methods of modern mathematical physics, v. IV, Acad. Press, NY, 1978.


\bibitem{Skrig_DAN2d} M.~M.~Skriganov, Proof of the Bethe-Sommerfeld conjecture in
dimension two, Soviet Math. Dokl. 20(1979), 956-959.

\bibitem{Skrig_book} M.~M.~Skriganov, Geometric and arithmetic methods in
the spectral theory of multidimensional periodic operators, Proc.
of the Steklov Inst. of Math. 171(1985), 1-117. Engl. transl. in
Proc. Steklov Inst. Math. 1987, no. 2.



\end{thebibliography}
\end{document}